\documentclass[amstex,12pt,russian,amssymb]{article}

\usepackage{mathtext}
\usepackage[cp1251]{inputenc}
\usepackage[T2A]{fontenc}
\usepackage[russian]{babel}
\usepackage[dvips]{graphicx}
\usepackage{amsmath}
\usepackage{amssymb}
\usepackage{amsxtra}
\usepackage{latexsym}
\usepackage{ifthen}

\begin{document}

\newcounter{lemma}
\newcommand{\lemma}{\par \refstepcounter{lemma}%
{\bf Лемма \arabic{lemma}.}}

\newcounter{corollary}
\newcommand{\corollary}{\par \refstepcounter{corollary}%
{\bf Следствие \arabic{corollary}.}}

\newcounter{remark}
\newcommand{\remark}{\par \refstepcounter{remark}%
{\bf Замечание \arabic{remark}.}}

\newcounter{theorem}
\newcommand{\theorem}{\par \refstepcounter{theorem}%
{\bf Теорема \arabic{theorem}.}}

\newcounter{proposition}
\newcommand{\proposition}{\par \refstepcounter{proposition}%
{\bf Предложение \arabic{proposition}.}}

\newcommand{\proof}{{\it Доказательство.\,\,}}
\renewcommand{\refname}{\centerline{\bf Список литературы}}

{\bf Е.А.~Петров} (Институт прикладной математики и механики НАН
Украины, г. Славянск)

{\bf Е.А.~Севостьянов} (Житомирский государственный университет
им.~И.~Франко)

\medskip\medskip
{\bf Є.О.~Петров} (Інститут прикладної математики і механіки НАН
України, м. Слов'янськ)

{\bf  Є.О.~Севостьянов} (Житомирський державний університет
ім.~І.~Франко)

\medskip\medskip
{\bf E.A.~Petrov} (Institute of Applied Mathematics and Mechanics of
NAS of Ukraine, Slavyansk)

{\bf E.A.~Sevost'yanov} (Zhitomir Ivan Franko State University)

\medskip
{\bf О равностепенной непрерывности классов Соболева в областях с
локально связной границей}

\medskip
Изучается поведение гомеоморфизмов классов Соболева в замыкании
заданной области. Доказано, что замыкании областей, границы которых
удовлетворяют определённым ограничениям, указанные классы
равностепенно непрерывны, как только их внутренняя дилатация порядка
$p$ имеет мажоранту конечного среднего колебания в каждой точке.

\medskip
{\bf Про одностайну неперервність класів Соболєва в областях з
локально зв'язною межею}

\medskip
Вивчається поведінка гомеоморфізмів класів Орліча--Соболєва в
замиканні заданої області. Доведено, що в замиканні областей, межі
яких задовольняють деякі обмеження, зазначені класи є одностайно
неперервними, як тільки їх внутряшня порядку $p$ має скінченне
середнє коливання в кожній точці.

\medskip
{\bf On equicontinuity of Sobolev classes in domains with locally
connected boundary}

\medskip
A behavior of homeomorphisms of Orlicz classes in a closure of a
domain is investigated. It is proved that above classes are
equicontinuous in the closure of domains with some restrictions on
it's boundaries provided that the corresponding inner dilatation of
order $p$ has a majorant of finite mean oscillation at every point.

\newpage
{\bf 1. Введение.} Настоящая работа посвящена изучению
гомеоморфизмов с конечным искажением, в частности,
$Q$-гомеоморфизмов, предложенных к изучению О. Мартио и изученных им
в совместной монографии \cite{MRSY} (см. также \cite{IM} и
\cite{GRSY}). Основные определения и обозначения, использующиеся
ниже, могут быть найдены в монографии \cite{MRSY} либо статье
\cite{Sev$_1$}.

В сравнительно недавней публикации одного из авторов Е.А.
Севостьянова установлена равностепенная непрерывность в замыкании
области некоторого класса гомеоморфизмов (см. \cite{Sev$_1$}). В
данной работе мы несколько усилим упомянутые результаты,
рассматривая более широкие классы отображений, при которых они всё
ещё имеют место. В качестве основных утверждений здесь предлагаются
результаты, относящиеся к классам Соболева $W_{loc}^{1, p},$
$p\geqslant 1,$ в то время как классы, аналогичные рассмотренным в
\cite{Sev$_1$}, приводятся в качестве вспомогательных.

Напомним некоторые определения, а также приведём формулировку
основных результатов работы. Всюду далее $D$ -- область в ${\Bbb
R}^n,$ $n\ge 2,$ $m$ -- мера Лебега в ${\Bbb R}^n.$ Пусть $U$ --
открытое множество, $U\subset {\Bbb R}^n,$ $u:U\rightarrow {\Bbb R}$
-- некоторая функция, $u\in L_{loc}^{\,1}(U).$ Предположим, что
найдётся функция $v\in L_{loc}^{\,1}(U),$ такая что
$$\int\limits_U \frac{\partial \varphi}{\partial x_i}(x)u(x)dm(x)=
-\int\limits_U \varphi(x)v(x)dm(x)$$
для любой функции $\varphi\in C_1^{\,0}(U).$ Тогда будем говорить,
что функция $v$ является {\it обобщённой производной первого порядка
функции $u$ по переменной $x_i$} и обозначать символом:
$\frac{\partial u}{\partial x_i}(x):=v.$ Функция $u\in
W_{loc}^{1,p}(U),$ $p\geqslant 1,$ если $u$ имеет обобщённые
производные первого порядка по каждой из переменных в $U,$ которые
являются локально интегрируемыми в $U$ в степени $p.$ Пусть $G$ --
открытое множество в ${\Bbb R}^n.$ Отображение $f:D\rightarrow {\Bbb
R}^n$ принадлежит {\it классу Соболева} $W^{1,p}_{loc}(G),$ пишем
$f\in W^{1,p}_{loc}(G),$ $p\geqslant 1,$ если все координатные
функции $f=(f_1,\ldots,f_n)$ обладают обобщёнными частными
производными первого порядка, которые локально интегрируемы в $G$ в
степени $p.$

Отображение $f:D\rightarrow \overline{{\Bbb R}^n}$ обладает {\it
$N$-свой\-с\-т\-вом} (Лузина), если из условия $m(E)=0$ следует, что
$m(f(E))=0.$ Отображение $f:D\rightarrow \overline{{\Bbb R}^n}$
обладает {\it $N^{\,-1}$-свойством}, если из условия $m(E)=0$
следует, что $m\left(f^{\,-1}(E)\right)=0.$

\medskip
Для отображений класса $W_{loc}^{1,1}$ и произвольного $p\geqslant
1$ корректно определена так называемая {\it внутренняя дилатация
$K_{I, p}(x,f)$ отображения $f$ порядка $p$ в точке $x$},
определяемая равенствами
\begin{equation}\label{eq0.1.1A}
K_{I, p}(x,f)\quad =\quad\left\{
\begin{array}{rr}
\frac{|J(x,f)|}{{l\left(f^{\,\prime}(x)\right)}^p}, & J(x,f)\ne 0,\\
1,  &  f^{\,\prime}(x)=0, \\
\infty, & \text{в\,\,остальных\,\,случаях}
\end{array}
\right.\,.
\end{equation}

Область $D$ называется {\it локально связной в точ\-ке}
$x_0\in\partial D,$ если для любой окрестности $U$ точки $x_0$
найдется окрестность $V\subset U$ точки $x_0$ такая, что $V\cap D$
связно, см. \cite[6.49.I]{Ku}. Следуя \cite[раздел 7.22]{He} будем
говорить, что борелева функция $\rho\colon  X\rightarrow [0,
\infty]$ является {\it верхним градиентом} функции $u\colon
X\rightarrow {\Bbb R},$ если для всех спрямляемых кривых $\gamma,$
соединяющих точки $x$ и $y\in X,$ выполняется неравенство
$|u(x)-u(y)|\leqslant \int\limits_{\gamma}\rho\,|dx|,$ где, как
обычно, $\int\limits_{\gamma}\rho\,|dx|$ обозначает линейный
интеграл от функции $\rho$ по кривой $\gamma.$ Будем также говорить,
что в указанном пространстве $X$ выполняется {\it $(1;
p)$-неравенство Пуанкаре,} если найдётся постоянная $C\geqslant 1$
такая, что для каждого шара $B\subset X,$ произвольной ограниченной
непрерывной функции $u\colon X\rightarrow {\Bbb R}$ и любого её
верхнего градиента $\rho$ выполняется следующее неравенство:
$$\frac{1}{\mu(B)}\int\limits_{B}|u-u_B|d\mu(x)\leqslant C\cdot({\rm diam\,}B)\left(\frac{1}{\mu(B)}
\int\limits_{B}\rho^p d\mu(x)\right)^{1/p}\,,$$
где $u_B:=\frac{1}{\mu(B)}\int\limits_Bu(x)d\mu(x).$ Метрическое
пространство $(X, d, \mu)$ назовём {\it $n$-регулярным по Альфорсу,}
если при каждом $x_0\in X,$ некоторой постоянной $C\geqslant 1$ и
всех $R<{\rm diam}\,X$
$$\frac{1}{C}R^{n}\leqslant \mu(B(x_0, R))\leqslant CR^{n}\,.$$

\begin{remark}\label{rem1}
Одним из примеров $n$-регулярного по Альфорсу пространства
относительно евклидовой метрики и меры Лебега в ${\Bbb R}^n,$ в
котором выполнено $(1; p)$-неравенство Пуанкаре, является единичный
шар ${\Bbb B}^n,$ см. \cite[предложение~2.1]{Sev$_6$}.
\end{remark}

\medskip Согласно \cite[разд.~6.1, гл.~6]{MRSY}, будем говорить, что
функция ${\varphi}:D\rightarrow{\Bbb R}$ имеет {\it конечное среднее
колебание} в точке $x_0\in D$, пишем $\varphi\in FMO(x_0),$ если
$\overline{\lim\limits_{\varepsilon\rightarrow
0}}\quad\frac{1}{\Omega_n\varepsilon^n} \int\limits_{B(x_0,
\varepsilon)}
|{\varphi}(x)-\overline{{\varphi}}_{\varepsilon}|\,dm(x)<\infty,
$
где $\Omega_n$ -- объём единичного шара ${\Bbb B}^n$ в ${\Bbb R}^n$
и
$\overline{{\varphi}}_{\varepsilon}\,=\,\frac{1}{\Omega_n\varepsilon^n}\int\limits_{B(
x_0,\,\varepsilon)} {\varphi}(x)\,dm(x).$

Пусть $(X,d)$ и $\left(X^{\,{\prime}},{d}^{\,{\prime}}\right)$ ~---
метрические пространства с расстояниями $d$  и ${d}^{\,{\prime}},$
соответственно. Семейство $\frak{F}$ отображений $f:X\rightarrow
{X}^{\,\prime}$ называется {\it равностепенно непрерывным в точке}
$x_0 \in X,$ если для любого $\varepsilon > 0$ найдётся $\delta
> 0,$ такое, что ${d}^{\,\prime} \left(f(x),f(x_0)\right)<\varepsilon$ для
всех $f \in \frak{F}$ и  для всех $x\in X$ таких, что
$d(x,x_0)<\delta.$ Говорят, что $\frak{F}$ {\it равностепенно
непрерывно}, если $\frak{F}$ равностепенно непрерывно в каждой точке
из $x_0\in X.$ Всюду далее, если не оговорено противное, $d$ -- одна
из метрик в пространстве простых концов относительно области $D,$
упомянутых выше, а $d^{\,\prime}$ -- евклидова метрика.

\medskip
Для числа $p\geqslant 1,$ областей $D,$ $D^{\,\prime}\subset {\Bbb
R}^n,$ $b_0\in D,$ $b_0^{\,\prime}\in D^{\prime}$ и произвольной
измеримой по Лебегу функции $Q(x): {\Bbb R}^n\rightarrow [0,
\infty],$ такой, что $Q(x)\equiv 0$ при $x\not\in D,$ обозначим
символом $\frak{F}_{p, b_0, b_0^{\,\prime}, Q}(D, D^{\,\prime})$
семейство всех гомеоморфизмов $f:D\rightarrow D^{\,\prime}$ класса
$W_{loc}^{1, p}$ в $D,$ $f(D)=D^{\,\prime},$ таких что $K_{I, p}(x,
f)\leqslant Q(x)$ и $f(b_0)=b_0^{\,\prime}.$ Справедливо следующее
утверждение.

\medskip
\begin{theorem}\label{th7}{\sl\,
Пусть $n\geqslant 2,$ $n-1<p\leqslant n,$ область $D$ имеет не менее
одной конечной граничной точки, $D$ локально связна на границе,
область $D^{\,\prime}$ ограничена и кроме того, $D^{\,\prime}$
является $n$-регулярным по Альфорсу пространством с евклидовой
метрикой и мерой Лебега в ${\Bbb R}^n$, в котором выполнено $(1;
p)$-неравенство Пуанкаре. Пусть также $Q\in FMO(\overline{D}),$ либо
$Q\in L_{loc}^1(D)$ и в каждой точке $x_0\in \overline{D}$ при
некотором $\varepsilon_0=\varepsilon_0(x_0)>0$ и всех
$0<\varepsilon<\varepsilon_0$
$$
\int\limits_{\varepsilon}^{\varepsilon_0}
\frac{dt}{t^{\frac{n-1}{n-p}}q_{x_0}^{\,\frac{1}{p-1}}(t)}<\infty\,,\qquad
\int\limits_{0}^{\varepsilon_0}
\frac{dt}{t^{\frac{n-1}{n-p}}q_{x_0}^{\,\frac{1}{p-1}}(t)}=\infty\,,
$$
где
$q_{x_0}(r):=\frac{1}{\omega_{n-1}r^{n-1}}\int\limits_{|x-x_0|=r}Q(x)\,d{\mathcal
H}^{n-1}.$
Тогда каждый элемент $f\in \frak{F}_{p, b_0, b_0^{\,\prime}, Q}(D,
D^{\,\prime})$ продолжается до непрерывного отображения
$\overline{f}:\overline{D}\rightarrow\overline{D^{\,\prime}}$, при
этом, семейство отображений $\overline{\frak{F}}_{p, b_0,
b_0^{\,\prime}, Q}(D, D^{\,\prime}),$ состоящее из всех продолженных
таким образом гомеоморфизмов, является равностепенно непрерывным, а
значит, и нормальным  в $\overline{D}$.}
\end{theorem}

\medskip
{\bf 2. Основные леммы. Равностепенная непрерывность кольцевых
гомеоморфизмов.} Пусть $E,$ $F\subset \overline{{\Bbb R}^n}$ --
произвольные множества. Обозначим через $\Gamma(E,F,D)$ семейство
всех кривых $\gamma:[a,b]\rightarrow\overline{{\Bbb R}^n},$ которые
соединяют $E$ и $F$ в $D,$ т.е. $\gamma(a)\in E,$ $\gamma(b)\in F$ и
$\gamma(t)\in D$ при $t\in (a, b).$ Дальнейшее изложение существенно
опираются на аппарат нак называемых $Q$-гомеоморфизмов относительно
$p$-модуля (см. \cite[глава~7]{MRSY}). Дадим определение этого
класса отображений. Пусть $x_0\in \overline{D},$ тогда отображение
$f:D\rightarrow \overline{{\Bbb R}^n}$ будем называть {\it кольцевым
$Q$-отображением относительно $p$-модуля в точке $x_0$,} если для
каждых $0<r_1<r_2<d_0:=\sup\limits_{x\in D}|x-x_0|$
\begin{equation}\label{eq14} M_p(f(\Gamma(S_1, S_2, D)))\leqslant \int\limits_{A(x_0,
r_1, r_2)\cap D}Q(x) \eta^p(|x-x_0|)dm(x)\,,
\end{equation}
где
\begin{equation}\label{eq1B}
A(x_0, r_1, r_2)=\{x\in {\Bbb R}^n: r_1<|x-x_0|<r_2\}\,,
\end{equation}
$S_i=S(x_0, r_i),$ $i=1,2,$ и $\eta: (r_1, r_2)\rightarrow
[0,\infty]$ -- произвольная измеримая по Лебегу функция такая, что
\begin{equation}\label{eq28*}
\int\limits_{r_1}^{r_2}\eta(r)dr=1\,.
\end{equation}

\medskip
Справедливо следующее утверждение (см.~\cite[предложение~4.7]{AS}).

\medskip
\begin{proposition}\label{pr2}
{\sl Пусть $X$ --- $\alpha$-регулярное по Альфорсу метрическое
пространство с мерой, в котором выполняется $(1; p)$-неравенство
Пуанкаре, $\alpha-1\leqslant p\leqslant \alpha.$ Тогда для
произвольных континуумов $E$ и $F,$ содержащихся в шаре $B(x_0, R),$
и некоторой постоянной $C>0$ выполняется неравенство
$M_p(\Gamma(E, F, X))\geqslant \frac{1}{C}\cdot\frac{\min\{{\rm
diam}\,E, {\rm diam}\,F\}}{R^{1+p-\alpha}}.$ }
\end{proposition}

\medskip
Имеет место следующее утверждение, обобщающее
\cite[лемма~3.1]{Sev$_1$} на случай произвольного порядка $p$ модуля
семейств кривых.

\medskip
\begin{lemma}\label{lem3}
{\sl\, Пусть $n\geqslant 2,$ $p\in (n-1, n],$ область $D\subset
{\Bbb R}^n$ локально связна в точках границы, $x_0\in
\partial D,$ $x_0\ne \infty,$ измеримая по
Лебегу функция $Q:{\Bbb R}^n\rightarrow [0, \infty]$ равна нулю вне
области $D,$ а отображение $f:D\rightarrow \overline{{\Bbb R}^n}$
является кольцевым $Q$-гомеоморфизмом в точке $x_0$ относительно
$p$-модуля, таким что $b_0^{\,\prime}=f(b_0)$ для некоторых $b_0\in
D$ и $b_0^{\,\prime}\in D^{\,\prime}=f(D).$ Пусть также
$D^{\,\prime}$ ограничена и является является $n$-регулярным по
Альфорсу пространством с евклидовой метрикой и мерой Лебега в ${\Bbb
R}^n$, в котором выполнено $(1; p)$-неравенство Пуанкаре.
Предположим, найдётся $\varepsilon_0=\varepsilon(x_0)>0,$
такое, что при некотором $0<p^{\,\prime}<p,$ $\varepsilon\rightarrow
0$ и некоторой постоянной $K>0$ выполнено условие
\begin{equation}\label{eq5***}
\int\limits_{A(x_0, \varepsilon, \varepsilon_0)}
Q(x)\cdot\psi^{\,p}(|x-x_0|)\,dm(x) \leqslant K\cdot
I^{p^{\,\prime}}(\varepsilon, \varepsilon_0)\,,
\end{equation}
где сферическое кольцо $A(x_0, \varepsilon, \varepsilon_0)$
определено как в (\ref{eq1B}), а $\psi$ -- некоторая заданная
неотрицательная измеримая функция, такая, что при всех
$\varepsilon\in(0, \varepsilon_0)$
\begin{equation}\label{eq7***}
I(\varepsilon,
\varepsilon_0):=\int\limits_{\varepsilon}^{\varepsilon_0}\psi(t)\,dt
< \infty\,,
\end{equation}
при этом, $I(\varepsilon, \varepsilon_0)\rightarrow \infty$ при
$\varepsilon\rightarrow 0.$ Тогда найдётся число
$\widetilde{\varepsilon_0}=\widetilde{\varepsilon_0}(x_0)\in (0,
\varepsilon_0)$ такое, что при каждом $\sigma\in (0,
\widetilde{\varepsilon_0})$ и любого континуума $E_1\subset B(x_0,
\sigma)\cap D$ выполнено неравенство
\begin{equation}\label{eq3.10}
{\rm diam\,}f(E_1)\leqslant C\cdot R^{1+p-n}\cdot K\cdot
I^{p^{\,\prime}-p}(\sigma, \varepsilon_0)\cdot\Delta(\sigma,
\widetilde{\varepsilon_0}, \varepsilon_0)\,,
\end{equation}
где ${\rm diam\,}f(E_1)$ -- евклидов диаметр множества $f(E_1),$
\begin{equation}\label{eq1.3}
\Delta(\sigma, \widetilde{\varepsilon_0}, \varepsilon_0)=\left(
1+\frac{\int\limits_{\widetilde{\varepsilon_0}}^{\varepsilon_0}\psi(t)\,dt}
{\int\limits_{\sigma}^{\widetilde{\varepsilon_0}}\psi(t)\,dt}\right)^{p}\,,\end{equation}
$\delta = \frac{1}{2}\cdot d\left(b_0^{\,\prime},
\partial D^{\,\prime}\right),$ $d$ -- евклидово расстояние
между множествами, $R$ -- радиус шара, содержащего область
$D^{\,\prime},$ а $C$ -- постоянная из предложения \ref{pr2}.}
\end{lemma}

\medskip
\begin{proof} Если граница
$\partial D$ области $D$ содержит ещё хотя бы одну точку $y_0\in
\partial D,$ $y_0\ne \infty,$ введём её в рассмотрение. В противном случае полагаем
$y_0:=\infty.$ Выберем $\varepsilon_1\in (0, \varepsilon_0)$ так,
чтобы точка $b_0$ может быть соединена непрерывной кривой с точкой
$y_0,$ лежащей целиком в $D\setminus B(x_0, \varepsilon_1),$ кроме
своей концевой точки $y_0.$ Заметим, что это возможно сделать, так
как область $D,$ по условию, предполагалась локально связной в любой
точке границы (см. \cite[предложение 13.2, гл.~13]{MRSY}). Последнюю
кривую обозначим через $E_3.$ Не ограничивая общности, можно
считать, что $I(\varepsilon_1, \varepsilon_0)>0,$ так как по условию
$I(\varepsilon, \varepsilon_0)\rightarrow \infty$ при
$\varepsilon\rightarrow 0.$ Хорошо известно, что предельное
множество $C(\partial G, g)$ лежит на границе области
$G^{\,\prime}:=g(G),$ как только $g:G\rightarrow G^{\,\prime}$ ~---
гомеоморфизм (см., напр., \cite[предложение~13.5, гл.~13]{MRSY}).
Учитывая сказанное выше, $C(y_0, f)\in
\partial D^{\,\prime}.$ Отсюда вытекает, что найдётся точка $a_0\in E_3,$
такая, что
\begin{equation}\label{eq46*}
d\left(b_0^{\,\prime}, f(a_0)\right)\ge \frac{1}{2}\cdot d
\left(b_0^{\,\prime}, \partial D^{\,\prime}\right):=\delta\,,
\end{equation}
где $d$ обозначает евклидово расстояние между множествами. Замкнутую
подкривую кривой $E_3,$ соединяющую точки $a_0$ и $b_0$ в $D,$
обозначим через $E_2.$ Пусть $\sigma\in (0, \varepsilon_1),$
$E_1\subset \overline{B(x_0, \sigma)}\cap D$ ~--- произвольный
континуум.

Рассмотрим измеримую функцию
$$\eta_{\sigma}(t)= \left\{
\begin{array}{rr}
\psi(t)/I(\sigma, \varepsilon_1), &   t\in (\sigma, \varepsilon_1),\\
0,  &  t\not\in (\sigma, \varepsilon_1)\,,
\end{array}
\right.$$
%
%
где, как и прежде, величина $I(a, b)$ определяется соотношением
$I(a, b)=\int\limits_a^b\psi(t)\,dt.$ Заметим, что функция
$\eta_{\sigma}(t)$ удовлетворяет соотношению вида (\ref{eq28*}), где
вместо $r_1$ и $r_2$ участвуют $\sigma$ и $\varepsilon_1,$
соответственно. Заметим, что $\Gamma\left(E_1, E_2,
D\right)>\Gamma(S(x_0, \sigma), S(x_0, \varepsilon_1), D)$ и значит,
$f(\Gamma\left(E_1, E_2, D\right))>f(\Gamma(S(x_0, \sigma), S(x_0,
\varepsilon_1),D)),$ откуда
$$M_p(f(\Gamma\left(E_1, E_2, D\right)))\leqslant M_p(f(\Gamma(S(x_0,
\sigma), S(x_0, \varepsilon_1), D))$$
(см. \cite[теорема~6.4]{Va}). В таком случае, ввиду условий
(\ref{eq5***})--(\ref{eq7***})
$$M_p\left(\Gamma\left(f(E_1), f(E_2),
D^{\,\prime}\right)\right)=$$
\begin{equation}\label{eq37***}
=M_p\left(f\left(\Gamma\left(E_1, E_2,
D\right)\right)\right)\leqslant M_p(f(\Gamma(S(x_0, \sigma), S(x_0,
\varepsilon_1), D )))\leqslant
\end{equation}
$$\leqslant\,\frac{K\cdot I^{p^{\,\prime}}(\sigma, \varepsilon_0)}{I^{p}(\sigma,
\varepsilon_1)}=K\cdot I^{p^{\,\prime}-p}(\sigma,
\varepsilon_0)\cdot\Delta(\sigma, \varepsilon_1, \varepsilon_0)\,,$$
где $\Delta$ определяется из соотношения (\ref{eq1.3}).
С другой стороны, по предложению \ref{pr2}
\begin{equation}\label{eq11}
\min\{{\rm diam\,}f(E_1), {\rm diam\,}f(E_2)\}\leqslant
CR^{1+p-n}M_{p}\left(\Gamma\left(f(E_1), f(E_2),
D^{\,\prime}\right)\right)\,,
\end{equation}
где $R$ -- радиус шара, содержащего область $D^{\,\prime}.$

Поскольку из (\ref{eq37***}) вытекает, что
$M_{p}\left(\Gamma\left(f(E_1), f(E_2),
D^{\,\prime}\right)\right)\rightarrow 0$ при $\sigma\rightarrow 0,$
то ввиду (\ref{eq11}) можно выбрать $\widetilde{\varepsilon_0}\in
(0, \varepsilon_1)$ таким, что что при каждом $\sigma\in (0,
\widetilde{\varepsilon_0})$ и любого континуума $E_1\subset B(x_0,
\sigma)\cap D$ будет выполнено следующее условие: $\min\{{\rm
diam\,}f(E_1), {\rm diam\,}f(E_2)\}<\delta.$ Но из (\ref{eq46*})
следует, что ${\rm diam\,}f(E_2)\geqslant \delta$ и, значит,
$\min\{{\rm diam\,}f(E_1), {\rm diam\,}f(E_2)\}={\rm diam\,}f(E_1).$
Тогда по неравенству (\ref{eq11})
\begin{equation}\label{eq11A}
{\rm diam\,}f(E_1)\leqslant CR^{1+p-n}M_{p}\left(\Gamma\left(f(E_1),
f(E_2), D^{\,\prime}\right)\right)\,.
\end{equation}
В этом случае, из (\ref{eq37***}) и (\ref{eq11A}) вытекает, что
\begin{equation}\label{eq13}
{\rm diam\,}f(E_1)\leqslant C\cdot R^{1+p-n}\cdot K\cdot
I^{p^{\,\prime}-p}(\sigma, \varepsilon_0)\cdot\Delta(\sigma,
\varepsilon_1, \varepsilon_0)\,.
\end{equation}
Осталось лишь заметить, что $\Delta(\sigma, \varepsilon_1,
\varepsilon_0)<\Delta(\sigma, \widetilde{\varepsilon_0},
\varepsilon_0),$ что следует прямо из определения величины $\Delta$
в (\ref{eq1.3}), так что из (\ref{eq13}) вытекает соотношение
(\ref{eq3.10}). Лемма доказана.~$\Box$
\end{proof}

\medskip
Для заданных областей $D,$ $D^{\,\prime}\subset {\Bbb R}^n,$ $n\ge
2,$ $n-1<p\leqslant n,$ измеримой по Лебегу функции $Q(x):{\Bbb
R}^n\rightarrow [0, \infty],$ $Q(x)\equiv 0$ при $x\not\in D,$
$b_0\in D,$ $b_0^{\,\prime}\in D^{\,\prime},$ обозначим через
$\frak{G}_{p, b_0, b_0^{\,\prime}, Q}\left(D, D^{\,\prime}\right)$
семейство всех кольцевых $Q$-гомеоморфизмов $f:D\rightarrow
D^{\,\prime}$ в $\overline{D}$ относительно $p$-модуля, таких что
$f(D)=D^{\,\prime},$ $b_0^{\,\prime}=f(b_0).$ В наиболее общей
ситуации основное утверждение настоящего раздела может быть
сформулировано следующим образом.

\medskip
\begin{lemma}\label{lem3A}
{\sl\, Пусть $p\in (n-1, n],$ область $D$ локально связна на границе
и имеет не менее одной конечной граничной точки, а область
$D^{\,\prime}$ ограничена, имеет локально квазиконформную границу и,
одновременно, является пространством $n$-регулярным по Альфорсу
относительно евклидовой метрики и меры Лебега в ${\Bbb R}^n,$ в
котором выполнено $(1; p)$-неравенство Пуанкаре.

Предположим, что для каждого $x_0\in \overline{D}$ найдётся
$\varepsilon_0=\varepsilon(x_0)>0,$
такое, что при некотором $0<p^{\,\prime}<p$ и
$\varepsilon\rightarrow 0$ выполнено условие (\ref{eq5***}), где
сферическое кольцо $A(x_0, \varepsilon, \varepsilon_0)$ определено
как в (\ref{eq1B}), а $\psi$ -- некоторая неотрицательная измеримая
функция, такая, что при всех $\varepsilon\in(0, \varepsilon_0)$
выполнено условие (\ref{eq7***}),
при этом, $I(\varepsilon, \varepsilon_0)\rightarrow \infty$ при
$\varepsilon\rightarrow 0.$

Тогда каждое $f\in\frak{G}_{p, b_0, b_0^{\,\prime}, Q}\left(D,
D^{\,\prime}\right)$ продолжается до непрерывного отображения
$f:\overline{D}\rightarrow \overline{D^{\,\prime}},$ при этом
семейство таким образом продолженных отображений является
равностепенно непрерывным в $\overline{D}.$ }
\end{lemma}

\begin{proof} Каждое отображение $f$ имеет непрерывное продолжение на
$\overline{D}$ в силу \cite[лемма~1]{Sev$_2$}. Равностепенная
непрерывность семейства $\frak{G}_{p, b_0, b_0^{\,\prime},
Q}\left(D, D^{\,\prime}\right)$ во внутренних точках области $D$
следует, например, из \cite[лемма~3]{Sev$_3$}.

Осталось показать равностепенную непрерывность  $\frak{G}_{p, b_0,
b_0^{\,\prime}, Q}\left(D, D^{\,\prime}\right)$ на $\partial D.$

Предположим противное, а именно, что семейство отображений
$\frak{G}_{p, b_0, b_0^{\,\prime}, Q}\left(D, D^{\,\prime}\right)$
не является равностепенно непрерывным в некоторой точке $x_0\in
\partial D.$ Тогда найдутся число $a>0,$ последовательность $y_k\in
\overline{D},$ $k=1,2,\ldots$ и элементы $f_k\in\frak{G}_{p, b_0,
b_0^{\,\prime}, Q}\left(D, D^{\,\prime}\right)$ такие, что $d(y_k,
x_0)<1/k$ и
\begin{equation}\label{eq6}
|f_k(y_k)-f_k(x_0)|\geqslant a\quad\forall\quad k=1,2,\ldots,\,.
\end{equation}
Ввиду возможности непрерывного продолжения каждого $f_k$ на границу
$D$ в терминах простых концов, для всякого $k\in {\Bbb N}$ найдётся
элемент $x_k\in D$ такой, что $d(x_k, y_k)<1/k$ и
$|f_k(x_k)-f_k(y_k)|<1/k.$ Тогда из (\ref{eq6}) вытекает, что
\begin{equation}\label{eq7}
|f_k(x_k)-f_k(x_0)|\geqslant a/2\quad\forall\quad k=1,2,\ldots,\,.
\end{equation}
Аналогично, в силу непрерывного продолжения отображения $f_k$ в
$\overline{D}$ найдётся последовательность $x_k^{\,\prime}\in D,$
$x_k^{\,\prime}\rightarrow x_0$ при $k\rightarrow \infty$ такая, что
$|f_k(x_k^{\,\prime})-f_k(x_0)|<1/k$ при $k=1,2,\ldots\,.$ Тогда из
(\ref{eq7}) вытекает, что
\begin{equation}\label{eq8}
|f_k(x_k)-f_k(x_k^{\,\prime})|\geqslant a/4\quad\forall\quad
k=1,2,\ldots\,,\,.
\end{equation}
Поскольку по условию область $D$ является локально связной на
границе, то найдётся последовательность окрестностей $V_m$ точки
$x_0,$ $m=1,2,\ldots,$ лежащих в шаре $B(x_0, 2^{-m}),$ такая что
$W_m:=V_m\cap D$ есть связное множество при каждом $m\in {\Bbb N},$
а так как $x_k$ и $x_k^{\,\prime}\rightarrow x_0$ при
$k\rightarrow\infty,$ то найдётся подпоследовательность номеров
$m_k,$ $k=1,2,\ldots,$ такая что $x_{m_k}\in W_k$ и
$x^{\,\prime}_{m_k}\in W_k.$ Соединим точки $x_{m_k}$ и
$x^{\,\prime}_{m_k}$ кривой $C_k,$ лежащей в $W_k.$ Тогда при
достаточно больших $k$ число $2^{\,-k}$ меньше числа
$\widetilde{\varepsilon_0}$ из леммы \ref{lem3}, поэтому в этой
лемме можно положить $E_1:=C_k.$ В таком случае, из леммы \ref{lem3}
вытекает, что
$$|f_{m_k}(x_{m_k})-f_{m_k}(x^{\,\prime}_{m_k})|\leqslant {\rm
diam\,}f_{m_k}(C_k)\leqslant$$
\begin{equation}\label{eq12}
\leqslant C\cdot R^{1+p-n}\cdot K\cdot I^{p^{\,\prime}-p}(2^{-k},
\varepsilon_0)\cdot\Delta(2^{-k}, \widetilde{\varepsilon_0},
\varepsilon_0)\rightarrow 0\,, \quad k\rightarrow\infty\,,
\end{equation}
что противоречит (\ref{eq8}). Полученное противоречие указывает на
то, что исходное предположение об отсутствии равностепенной
непрерывности семейства $\frak{G}_{p, b_0, b_0^{\,\prime},
Q}\left(D, D^{\,\prime}\right)$ было неверным.~$\Box$
\end{proof}

\medskip
Основное утверждение настоящего раздела может быть сформулировано
следующим образом.

\medskip
\begin{theorem}\label{th1}
{\sl\, Пусть $p\in (n-1, n],$ область $D$ локально связна на границе
и имеет не менее одной конечной граничной точки, а область
$D^{\,\prime}$ ограничена, имеет локально квазиконформную границу и,
одновременно, является пространством $n$-регулярным по Альфорсу
относительно евклидовой метрики и меры Лебега в ${\Bbb R}^n,$ в
котором выполнено $(1; p)$-неравенство Пуанкаре.

Предположим,

1) $Q\in FMO(\overline{D}),$ либо

2) в каждой точке $x_0\in \overline{D}$ при некотором
$\varepsilon_0=\varepsilon_0(x_0)>0$ и всех
$0<\varepsilon<\varepsilon_0$
$$
\int\limits_{\varepsilon}^{\varepsilon_0}
\frac{dt}{t^{\frac{n-1}{n-p}}q_{x_0}^{\,\frac{1}{p-1}}(t)}<\infty\,,\qquad
\int\limits_{0}^{\varepsilon_0}
\frac{dt}{t^{\frac{n-1}{n-p}}q_{x_0}^{\,\frac{1}{p-1}}(t)}=\infty\,,
$$
где
$q_{x_0}(r):=\frac{1}{\omega_{n-1}r^{n-1}}\int\limits_{|x-x_0|=r}Q(x)\,d{\mathcal
H}^{n-1}.$

Тогда каждое $f\in\frak{G}_{p, b_0, b_0^{\,\prime}, Q}\left(D,
D^{\,\prime}\right)$ продолжается до непрерывного отображения
$f:\overline{D}\rightarrow \overline{D^{\,\prime}},$ при этом
семейство таким образом продолженных отображений является
равностепенно непрерывным в $\overline{D}.$ }
\end{theorem}

\medskip
{\it Доказательство теоремы \ref{th1}} сводится к тому, что условия
1) и 2), накладывающие ограничения на функцию $Q,$ влекут выполнение
условий (\ref{eq5***})--(\ref{eq7***}). Действительно, пусть вначале
имет место условие 1). Не ограничивая общности, можно считать, что
$x_0=0.$ Пусть $Q\in FMO(0)$ и $\varepsilon_0<\min\left\{\,\,{\rm
dist\,}\left(0, \,\partial D\right),\,\, e^{\,-1}\right\}.$ На
основании \cite[следствие 6.3, гл. 6]{MRSY} для функции
$0\,<\,\psi(t)\,=\,\frac
{1}{\left(t\,\log{\frac1t}\right)^{n/p}}$ будем иметь, что %
$$\int\limits_{\varepsilon<|x|<\varepsilon_0} Q(x)\cdot\psi^p(|x|)
 \ dm(x)\,= \int\limits_{\varepsilon<|x|< {\varepsilon_0}}\frac{Q(x)\,
dm(x)} {\left(|x| \log \frac{1}{|x|}\right)^n} = O \left(\log\log
\frac{1}{\varepsilon}\right)$$
%
при  $\varepsilon \rightarrow 0.$
Заметим также, что при указанных выше $\varepsilon$ выполнено
$\psi(t)\ge \frac {1}{t\,\log{\frac1t}},$ поэтому
$I(\varepsilon,
\varepsilon_0)\,:=\,\int\limits_{\varepsilon}^{\varepsilon_0}\psi(t)\,dt\,\ge
\log{\frac{\log{\frac{1}
{\varepsilon}}}{\log{\frac{1}{\varepsilon_0}}}}.$ Тогда
$$\int\limits_{\varepsilon<|x|<\varepsilon_0}
Q(x)\cdot\psi^p(|x|)
 \ dm(x)\le C\cdot I(\varepsilon, \varepsilon_0)$$ при $\varepsilon\rightarrow
 0.$ Таким образом, мы имеем соотношения вида
(\ref{eq5***})--(\ref{eq7***}), а, значит, в ситуации 1) заключение
теоремы \ref{th1} вытекает из леммы \ref{lem3A}.

Покажем, что аналогичное заключение можно сделать и в случае 2).
Положим в этом случае
\begin{equation}\label{eq2.3.9}
\psi(t)= \left \{\begin{array}{rr}
1/[t^{\frac{n-1}{p-1}}q_{x_0}^{\frac{1}{p-1}}(t)]\ , & \ t\in
(0,\varepsilon_0)\ ,
\\ 0\ ,  &  \ t\notin (0,\varepsilon_0)\ .
\end{array} \right.
\end{equation}
Тогда ввиду теоремы Фубини
\begin{equation*}
\int\limits_{\varepsilon<|x-x_0|<\varepsilon_0}Q(x)\cdot\psi^p(|x-x_0|)
\ dm(x)=\omega_{n-1}\,\int\limits_{\varepsilon}^{\varepsilon_0}\
\frac{dr}{r^{\frac{n-1}{p-1}}
q_{x_0}^{\frac{1}{p-1}}(r)}=\omega_{n-1}\cdot I(\varepsilon,
\varepsilon_0)\,,
\end{equation*}
где, как и прежде, величина $I(a, b)$ определяется соотношением
$I(a, b)=\int\limits_a^b\psi(t)\,dt,$ а $\omega_{n-1}$ -- площадь
единичной сферы в ${\Bbb R}^n.$ Следовательно, вновь выполнены
соотношения (\ref{eq5***})--(\ref{eq7***}), а, значит, в ситуации 2)
мы также получаем интересующее нас заключение из леммы
\ref{lem3A}.~$\Box$

\medskip
{\bf 4. Доказательство основного результата -- теоремы \ref{th7}.}
Для доказательства теоремы \ref{th7} покажем, что каждое отображение
$f\in\frak{F}_{p, b_0, b_0^{\,\prime}, Q}(D, D^{\,\prime})$
принадлежит также классу $\frak{G}_{p, b_0, b_0^{\,\prime},
Q}\left(D, D^{\,\prime}\right)$ (см. определения этих классов перед
теоремой \ref{th7} и перед леммой \ref{lem3A}). Действительно, пусть
$f\in\frak{F}_{p, b_0, b_0^{\,\prime}, Q}(D, D^{\,\prime}),$ тогда в
частности, $f\in W_{loc}^{1, p},$ $p>n-1.$ Отсюда ввиду
\cite[теорема 1.1]{Zi} $f^{\,-1}\in W_{loc}^{1, 1}.$ Поскольку
отображение $f$ обладает $N$-свойством Лузина, в силу замены
переменной под знаком интеграла (см. \cite[теорема~3.2.5]{Fe}) для
любого компакта $K\subset f(D)$
$$\int\limits_{K} \Vert g^{\,\prime}(y)\Vert^p
dm(y)=\int\limits_{f^{\,-1}(K)} K_{I, p}(x,
f)dm(x)<\int\limits_{f^{\,-1}(K)} Q(x)dm(x)<\infty\,,$$
где $g(y)=f^{\,-1}(y),$ откуда $f^{\,-1}\in W_{loc}^{1, p}.$ Кроме
того, так как $f\in W_{loc}^{1, p},$ $p>n-1,$ то $f$ является
дифференцируемым почти всюду (см. \cite[лемма~3]{Va$_2$}).

Таким образом, $f$ является дифференцируемым почти всюду
отображением, обладающим $N$- и $N^{\,-1}$-свойствами Лузина, таким
что $f^{\,-1}\in W_{loc}^{1, p}$ и $K_{I, p}(x, f)\leqslant Q(x),$ а
указанные отображения удовлетворяют неравенствам вида (\ref{eq14}),
что следует из \cite[теорема~2.2]{SalSev}. Таким образом,
$\frak{F}_{p, b_0, b_0^{\,\prime}, Q}(D, D^{\,\prime})\subset
\frak{G}_{p, b_0, b_0^{\,\prime}, Q}\left(D, D^{\,\prime}\right)$ и,
значит, заключение теоремы \ref{th7} вытекает из теоремы
\ref{th1}.~$\Box$

\medskip
КОНТАКТНАЯ ИНФОРМАЦИЯ

\medskip
\noindent{{\bf Евгений Александрович Петров}\\
Институт прикладной математики и механики
\\ НАН Украины, комн. 417\\
ул. генерала Батюка, 19 \\
г. Славянск, Украина, 84 116\\
e-mail: eugeniy.petrov@gmail.com

\medskip
\noindent{{\bf Евгений Александрович Севостьянов} \\
Житомирский государственный университет им.\ И.~Франко\\
ул. Большая Бердичевская, 40 \\
г.~Житомир, Украина, 10 008 \\ тел. +38 066 959 50 34 (моб.),
e-mail: esevostyanov2009@mail.ru}

\end{document}